\documentclass[10pt]{article}
\font\smallit=cmti10
\font\smalltt=cmtt10

\usepackage{amssymb,latexsym,amsmath,epsfig,amsthm}

\makeatletter

\renewcommand\section{\@startsection {section}{1}{\z@}
{-30pt \@plus -1ex \@minus -.2ex}
{2.3ex \@plus.2ex}
{\normalfont\normalsize\bfseries\boldmath}}

\renewcommand\subsection{\@startsection{subsection}{2}{\z@}
{-3.25ex\@plus -1ex \@minus -.2ex}
{1.5ex \@plus .2ex}
{\normalfont\normalsize\bfseries\boldmath}}

\renewcommand{\@seccntformat}[1]{\csname the#1\endcsname. }

\makeatother

\newtheorem*{dfn}{Definition}
\newtheorem{thm}{Theorem}
\newtheorem{lem}{Lemma}

\newtheorem{prop}{Proposition}
\newtheorem{cor}{Corollary}

\begin{document}

\begin{center}
\uppercase{\bf A Note on the Erd\H{o}s-Straus Conjecture}
\vskip 20pt
{\bf Kyle Bradford}
{\smallit Department of Mathematics and Statistics, Georgia Southern University, Statesboro, GA 30458, USA}\\
{\tt kbradford@georgiasouthern.edu}\\
\end{center}
\vskip 30pt
\centerline{\smallit Received: , Revised: , Accepted: , Published: } 
\vskip 30pt

\centerline{\bf Abstract}
\noindent This paper makes a fundamental assertion about the Erd\H{o}s-Straus  conjecture.  Suppose that for a prime $p$  there exists $x,y,z \in \mathbb{N}$  with $x \leq y \leq z$  so that $$ \frac{4}{p} = \frac{1}{x} + \frac{1}{y} + \frac{1}{z}. $$

\ 

The main contribution of this paper is that, under this assumption, the Erd\H{o}s-Straus conjecture can be reduced by one variable.  For example, it is necessarily true that $$ z = \frac{xyp}{\gcd(y,p) \gcd \left( xy, x+y \right)}.$$ 

\ 

Considering other reductions of the Erd\H{o}s-Straus conjecture, this paper suggests a method for proof.

\pagestyle{myheadings} 
\markright{\smalltt INTEGERS: 18 (2018)\hfill} 
\thispagestyle{empty} 
\baselineskip=12.875pt 
\vskip 30pt

\section{Preliminaries}

\noindent The ancient Egyptian society revered unit fractions as the fundamental building blocks of all other fractions.  In fact, ancient Egyptians ascribed special unit fractions to the eye of the god Horus.  Thus the name ``Egyptian fraction" was given to sums of unit fractions with positive denominators.  The monumental work of Paul Erd\H{o}s left a variety of unsolved problems related to Egyptian fractions \cite{rg}. This has led to a large body of work \cite{aaa, bh, ec, gm, rav, wweb2}.  One of the more common areas of study is the Erd\H{o}s-Straus and related conjectures \cite{bi, et, pe1, iw, dl, mo, ob, ro, san1, san2, san3, s, t, v, wweb1, wweb3, y, yang}.  In its original form it reduces to the following: given a prime number $p$  there exist natural numbers $x,y,$  and $z$  (w.l.o.g. $x \leq y \leq z$)  so that the Erd\H{o}s-Straus equation is solved.  This is expressed as
\begin{equation} \label{eqn: se}
\frac{4}{p} = \frac{1}{x} + \frac{1}{y} + \frac{1}{z}.
\end{equation}

The most famous approach is by Rosati \cite{ro} and is outlined in a book by Mordell \cite{mo}, which described modular identities for specific cases of primes, and an assumption that quadratic residues play a role in the solution of the problem.  More recently Elsholtz and Tao \cite{et}  attempted to find a solution using arithmetic number theory.  Both of these texts motivated a different approach in the paper presented here.  This theoretical paper is a continuation of the work in \cite{bi}, which was more computational in nature.  The patterns discovered here help simplify the problem by one dimension but do not solve the conjecture.

Assume that the conjecture is true.  Discerning a pattern between a given prime number $p$  and the associated solution values $x,y,$  and $z$  will determine the necessary conditions to show that a solution exists.  To prevent any ambiguity, note that this paper does not show that a solution exists.  Going forward, reserve $x,y,$  and $z$  as the general solution values for a prime $p$ and insist that $x \leq y \leq z$.  The following propositions are taken from an article \cite{et}.  


\begin{prop} \label{prop: a}
A prime number $p$ must divide at least one of its solution values $x,y,$  or $z$.
\end{prop}

\begin{prop} \label{prop: b}
For a given prime number $p$,  the solution values $x,y,$  and $z$  cannot simultaneously be divisible by $p$.
\end{prop}

The first lemma serves to further illuminate the nature of the solution values for a given prime $p$.

\begin{lem} \label{lem: a}
\noindent Given a prime number $p$, necessary conditions for the solution values $x \leq y \leq z$  are

\begin{equation} \label{reg: one}
\left\lceil \frac{p}{4} \right\rceil \leq x \leq \left\lfloor \frac{3p}{4} \right\rfloor \qquad  and  \qquad \left\lceil \frac{xp}{4x - p} \right\rceil \leq y \leq \left\lfloor \frac{2xp}{4x - p} \right\rfloor .
\end{equation}
\end{lem}

\ 

\noindent The necessary conditions (\ref{reg: one})  establish the subsequent lemmata. 

\begin{lem} \label{lem: b}
A prime number $p$  is coprime to the smallest solution value $x$.
\end{lem}

\begin{lem} \label{lem: c}
For a prime number $p$, both $\gcd (y,p^{2})\neq p^{2}$ and $\gcd (z,p^{2}) \neq p^{2}$.
\end{lem}

\noindent The fourth lemma guarantees that the solution value $z$  is divisible by $p$.

\begin{lem} \label{lem: d}
For a prime number $p$, if $\gcd (y,p) = p$, then $\gcd (z,p) = p$.   
\end{lem}

\begin{dfn} \label{def: zwei}
Assume that a solution to (\ref{eqn: se})  exists for a prime $p$.  Define a Type I solution as one so that $\gcd (x,p) = 1$, $\gcd (y,p) =1$,  and $\gcd (z,p) = p$.  Define a Type II solution as one so that $\gcd(x,p) = 1$, $\gcd (y,p) = p$,  and $\gcd (z,p) = p$.
\end{dfn}

This language is used in \cite{et} to describe types of solutions; however, the current paper insists on an ordering of the solution values.  Once classified into types, the solution values must be factored into their smallest relevant components. 

\begin{figure}[h] \label{fig: one}
\begin{centering}
\includegraphics[trim={0cm 12cm 3.2cm 1cm}, scale=0.30]{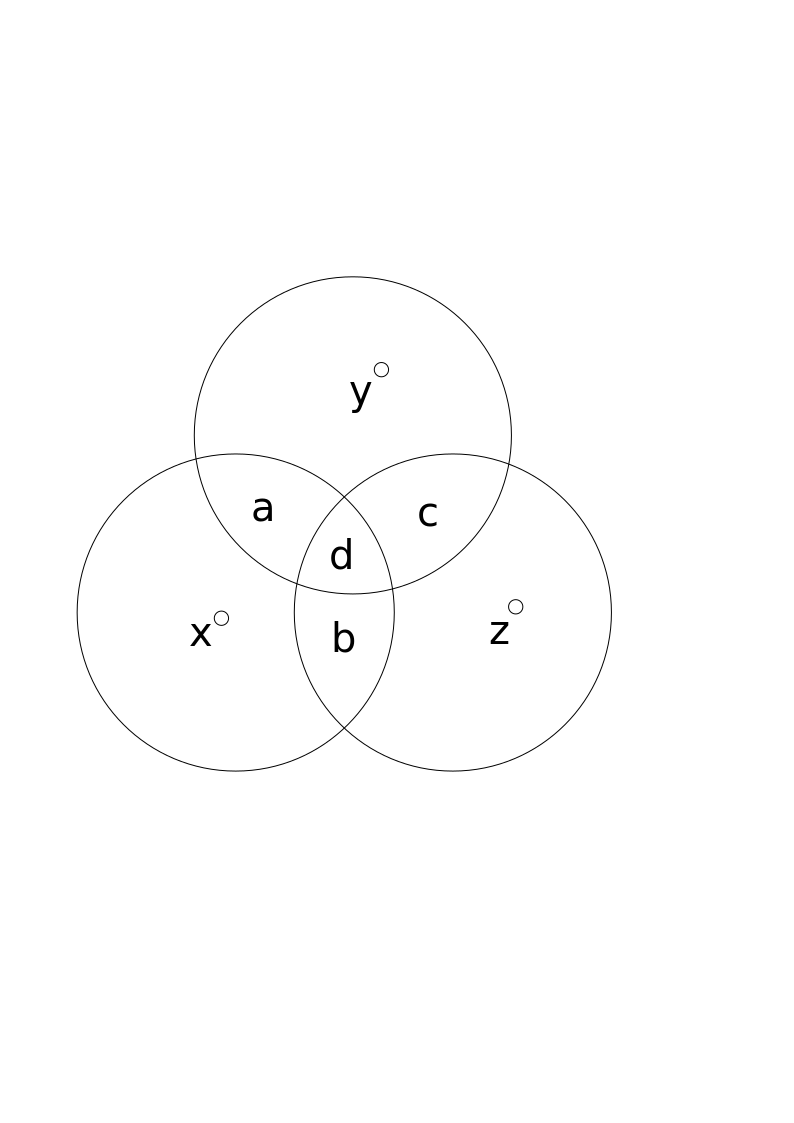}
\caption{This is a diagram of $x,y$, and $z$ into its factored parts.}
\end{centering}
\end{figure}

\begin{dfn} \label{def: drei}
Assume that a solution to (\ref{eqn: se})  exists for a prime $p$.  Reserve  $d = \gcd(x,y,z)$, $ a = \gcd(x,y) \slash d$, $b=\gcd(x,z) \slash d$,  and $c=\gcd(y,z) \slash d$.  Also reserve $x^{\circ}, y^{\circ}$,  and $z^{\circ}$  to be positive integers so that $x = x^{\circ} abd$, $y=y^{\circ}acd$,  and $z=z^{\circ}bcd$.
\end{dfn}

This makes clear that $a,b$,  and $c$  are pairwise coprime.  For Type I solutions $p|z^{\circ}$,  and for Type II solutions $p|c$.  Refer to Figure \ref{fig: one} for clarity. The final lemma reduces the complexity of the factorizations of $x,y$,  and $z$. 

\begin{lem} \label{lem: e}
For a prime number $p$  a Type I solution has $x^{\circ} = y^{\circ}=1$  and $z^{\circ} = p$, and a Type II solution has $x^{\circ} = y^{\circ} = z^{\circ} = 1$.
\end{lem}

These preliminary results assure that the factorizations of $x,y$,  and $z$  are fundamental for both types of solutions.  The results in this paper use these factorizations to reveal patterns that reduce the dimensionality of this conjecture.

\section{Results}

\noindent This section outlines the main results of this paper and provides motivation to a method of solution for the conjecture.

\begin{thm} \label{thm: a}
For a prime number $p$, if $x$  and $y$  are solution values to (\ref{eqn: se}), then necessarily

\begin{equation} \label{eqn: a}
4xy - (x+y)p = \gcd(y,p) \gcd ( xy, x+y ).
\end{equation}
\end{thm}

\noindent The most relevant implication from (\ref{eqn: a}) is that  

\begin{equation} \label{eqn: d}
z = \frac{xyp}{\gcd (y,p) \gcd ( xy, x+y )}.
\end{equation}

Notice $z$, by its definition in (\ref{eqn: d}), has to be an integer.  For a given $p$  it suffices to find $x,y \in \mathbb{N}$  in the region defined in (\ref{reg: one}) so that (\ref{eqn: a}) holds.  This, in essence, reduces the dimensionality of the problem by one variable.  This is also incredibly important because it reveals that the true nature of this problem depends on the gcd of the product of two numbers and the sum of those same two numbers.  This is not immediately understood from the original description, and this paper hopes to inform and motivate mathematicians who have studied this problem. The next theorem is motivated by the results in \cite{bi}.

\begin{figure}[h] \label{fig: two}
\fbox{\includegraphics[scale=0.70]{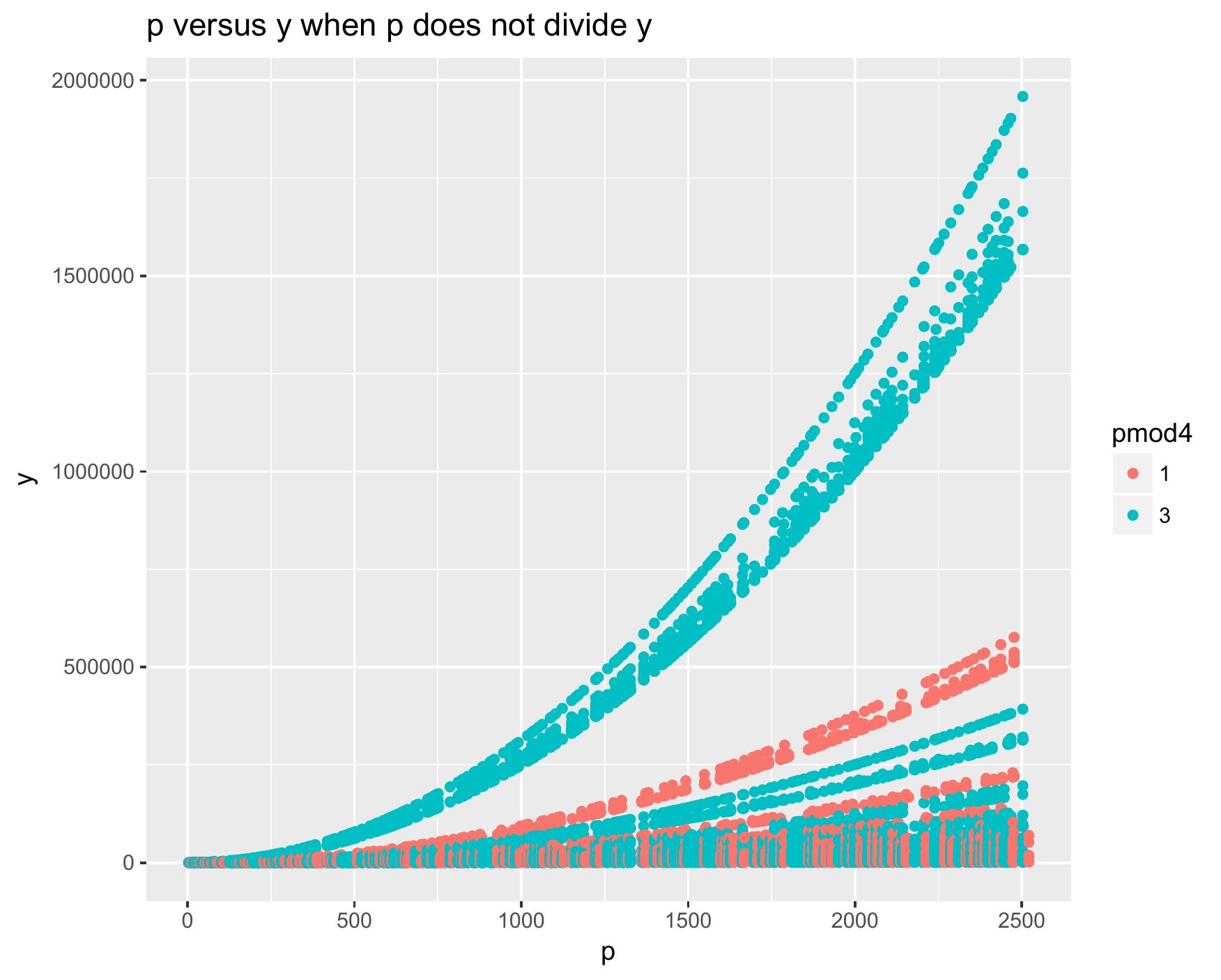}}
\caption{Type I solutions.  Red points denote that the prime has remainder 1 after dividing by 4, and blue points denote that the prime has remainder 3 after dividing by 4.}
\end{figure}

\begin{thm} \label{thm: b}
For a prime number $p$  with a Type I solution 
\begin{equation} \label{eqn: e} 
x = \left\lceil \frac{yp}{4y-p} \right\rceil .
\end{equation}
\end{thm}

This result is astounding.  The dimensionality of the problem is again reduced by one degree for a Type I solution.  It also holds for a vast majority of Type II solutions.  Computational evidence suggests that (\ref{eqn: e}) holds for over 97\%  of all solutions \cite{bi}.  Theorem \ref{thm: b}  also helps to reduce the bounds to Lemma \ref{lem: a}. \\

\begin{cor} \label{cor: a}
\noindent Given a prime number $p$, necessary conditions for the solution values $x \leq y \leq z$  are

\begin{equation} \label{reg: one}
\left\lceil \frac{p}{4} \right\rceil \leq x \leq \left\lceil \frac{p}{2} \right\rceil \qquad  and  \qquad \left\lceil \frac{xp}{4x - p} \right\rceil \leq y \leq \left\lfloor \frac{2xp}{4x - p} \right\rfloor .
\end{equation}
\end{cor}

Further note that for $3 \equiv p \mod 4$  that a Type I solution exists when $x=y=\left\lceil p \slash 2 \right\rceil$.  These are the only types of solutions that exist when $x = \left\lceil p \slash 2 \right\rceil$.  Considering the bounds in Corollary \ref{cor: a} and the results in Theorems \ref{thm: a}  and \ref{thm: b}, for a prime $p$ it suffices to find $y \in \mathbb{N}$  so that 

\begin{equation} \label{reg: three}
 \left\lceil \frac{p}{2} \right\rceil \leq y \leq \left\lfloor \frac{2 \left\lceil \frac{p}{4} \right\rceil p}{4 \left\lceil \frac{p}{4} \right\rceil - p} \right\rfloor 
\end{equation} 

\noindent and 
\begin{equation} \label{eqn: f}
(4y-p) \left\lceil \frac{yp}{4y-p} \right\rceil - yp = \gcd(y,p) \gcd \left(y \cdot \left\lceil \frac{yp}{4y-p} \right\rceil , y + \left\lceil \frac{yp}{4y-p} \right\rceil  \right).
\end{equation}

For $p$  large enough, it suffices to find a functional expression for $y$  that depends solely on $p$ so that (\ref{reg: three})  and (\ref{eqn: f}) hold.  It's important to note that the functional description for $y$  would have to lie between linear and quadratic behavior in $p$, although finding the correct description has proven elusive.   Figure \ref{fig: two} shows that many patterns exist between $y$  and $p$ for Type I solutions.  These patterns are found as modular identities outlined in previous papers \cite{mo}, but I maintain hope that a general pattern can be found.



%
%
%
%

\section*{Proofs}

\begin{proof} \label{lem: one}
\underline{Lemma 1}: \\

\noindent  First, assume that $x < p \slash 4$.  Substitute the bound $p \slash 4$  for $x$  in the algebraic expression $4xy - (x+y)p$  to see that $4xy - (x+y)p \leq - p^{2} \slash 4 < 0$.  Also notice that $xyp>0$  for all solution values.  By definition $z = xyp \slash (4xy - (x+y)p)$.  These observations imply that $z< 0$, but to be a solution value requires that $z>0$.  This creates a contradiction negating the assumption and implies that $x \geq p \slash 4$.  Recall that $x \in \mathbb{N}$ to conclude that $x \geq \lceil p \slash 4 \rceil$.  

Next, assume that $y > 2xp \slash (4x - p)$.  The previous result establishes that $4x - p > 0$.  When either side of the inequality is multiplied by $4x-p$,  the result is $4xy - yp > 2xp$.  This can be rewritten as $4xy - (x+y)p > xp$, which indicates that $4xy - (x+y)p > 0$ because $xp > 0$.  Multiplying either side of the previous inequality by $y>0$  implies that $(4xy - (x+y)p)y > xyp$.  Dividing both sides of the inequality by $4xy - (x+y)p$  shows that $y > xyp \slash (4xy - (x+y)p) = z$.  Solution values require that $y \leq z$.  This creates a contradiction that negates the assumption and implies that $y \leq 2xp \slash (4x-p)$.  Recall that $y \in \mathbb{N}$  to conclude that $y \leq \lfloor 2xp \slash (4x - p) \rfloor$.  

Next, assume that $y < xp \slash (4x - p)$.  Remember that $4x - p > 0$.  Multiplying either side of the inequality by  $4x-p$  enables the derivation of $4xy - (x+y)p < 0$.  Also remember that $xyp >0$.  Dividing $xyp$  by $4xy - (x+y)p$  shows that $z = xyp \slash (4xy - (x+y)p) < 0$.  To be a solution value requires that $z>0$.  This creates a contradiction that negates the assumption and implies that $y \geq xp \slash (4x - p)$.  Recall that $y \in \mathbb{N}$  to conclude that $y \geq \lceil xp \slash (4x - p) \rceil$. 

Finally, assume that $x > 3p \slash 4$.  For solution values the inequality $x\leq y \leq z$ makes both $y > 3p \slash 4$  and $z > 3p \slash 4$.  It implies that $4 \slash p > 1 \slash x + 1 \slash y + 1 \slash z$.  This contradicts the assumption that $x,y$,  and $z$  are solution values.  This contradiction negates the assumption and implies that $x \leq 3p \slash 4$.  Recall that $x \in \mathbb{N}$  to conclude that $x \leq \lfloor 3p \slash 4 \rfloor$. \\
\end{proof}


\begin{proof} \label{lem: two}
\underline{Lemma 2}: \\

\noindent Assume that $\gcd(x,p) \neq 1$.  This implies that $p|x$.  Lemma~\ref{lem: a}  shows that $x  \leq \lfloor 3p \slash 4 \rfloor$.  It is clear that $x < p$, so $p \nmid x$  which creates a contradiction.  This contradiction negates the assumption and implies that $x$ and $p$  are coprime.
\end{proof}


\begin{proof} \label{lem: three}
\underline{Lemma 3}: \\

\noindent First, assume that $\gcd ( y, p^{2}) = p^{2}$.  This implies that $p^{2} | y$.  Recall from Lemma~\ref{lem: a}  that $y \leq 2xp \slash (4x - p)$.  Start by showing that $2xp \slash (4x-p)$  decreases as the integer $x \geq \lceil p \slash 4 \rceil$  increases.  Consider that $2xp(4x-p) + 8xp > 2xp(4x-p)+ 8xp -2p^{2}$  because $0 > -p^{2}$, and rewrite it as $2xp(4(x+1)-p) > 2(x+1)p(4x-p)$.  Lemma~\ref{lem: a} also shows that $x > p \slash 4$.  As a consequence, both $4x-p>0$ and $4(x+1)-p>0$.  Dividing both sides of the previous inequality by $(4x-p)(4(x+1)-p)$  shows that $2xp \slash (4x-p) > 2(x+1)p \slash (4(x+1)-p)$.    This shows that $2xp \slash (4x-p)$  decreases as the integer $x \geq \lceil p \slash 4 \rceil$  increases.  Plugging the smallest possible integer value for $x$  into the expression $2xp \slash (4x-p)$  maximizes its possible value over integers $x \geq \lceil p \slash 4 \rceil$.  This implies that $y \leq 2 \lceil p \slash 4 \rceil p \slash (4 \lceil p \slash 4 \rceil -p)$.  Notice that $\lceil p \slash 4 \rceil \leq (p+1) \slash 4$ and $4 \lceil p \slash 4 \rceil - p \geq 1$  for all primes $p$.  This shows that $y \leq p(p+1) \slash 2 < p^{2}$.  If $y < p^{2}$, then $p^{2} \nmid y$.  which creates a contradiction.  This contradiction negates the assumption and implies that $\gcd (y, p^{2}) \neq p^{2}$. 

Next, assume that $\gcd(z,p^{2})=p^{2}$.  This implies that $p^{2} | z$.  Let $z^{*} \in \mathbb{N}$  so that $z=z^{*}p^{2}$.  Solution values necessarily satisfy the equation $y = xzp \slash (4xz - (x+z)p)$  which can be written as $y = xz^{*}p^{2} \slash ((4xz^{*} - z^{*}p)p - x)$.  Lemma~\ref{lem: b} shows that $\gcd(x,p)=1$  which implies that $p \nmid ((4xz^{*} - z^{*}p)p - x)$.  This implies that $p^{2} | y$, which further implies that $\gcd(y,p^{2}) = p^{2}$.  This creates a contradiction because it was just shown that $\gcd (y,p^{2}) \neq p^{2}$.  This contradiction negates the assumption and implies that $\gcd(z,p^{2}) \neq p^{2}$.
\end{proof}


\begin{proof} \label{lem: four}
\underline{Lemma 4}: \\

\noindent Let $\gcd (y,p) = p$  and for sake of contradiction assume that $\gcd (z,p) = 1$.  Let $y^{*} \in \mathbb{N}$  so that $y = y^{*}p$.  Solution values necessarily satisfy the equation $z = xyp \slash (4xy - (x+y)p)$  which can be written as $ z  =  xy^{*}p \slash ((4y^{*} - 1)x - y^{*}p)$.  Lemma~\ref{lem: b} shows that $\gcd(x,p)=1$  and Lemma~\ref{lem: c} shows that $\gcd(y^{*},p)=1$.  The assumption that $\gcd (z,p) =1$  requires that $p | (4y^{*} - 1)$.  Recalling the proof of Lemma~\ref{lem: c}, it was shown that $y \leq p(p+1) \slash 2$.  This implies that $y^{*} \leq (p+1) \slash 2$  which further implies that $4y^{*} - 1 \leq 2p + 1$.  If $p | (4y^{*} - 1)$  then either $4y^{*} - 1 = p$  or $4y^{*} - 1 = 2p$.  $4y^{*} - 1$  cannot be even, so $4y^{*} - 1 \neq 2p$.  If $4y^{*} - 1 = p$  then it is similarly clear that $3 \equiv p \mod 4$.  Writing $y^{*} = (p+1) \slash 4$  implies that $ z = (p+1)x \slash (4x-(p+1))$ which exists, is positive, and is maximized at the integer value $x = (p+5) \slash 4 \geq \left\lceil p \slash 4 \right \rceil$. This implies that $ z \leq ((p+1) \slash4) ((p+5) \slash 4)$.  Because $(p+5) \slash 4 < p$  for all primes such that $3 \equiv p \mod 4$,  this implies that $z < p(p+1) \slash 4$  which further implies that $z < y$.  This creates a contradiction to $x \leq y \leq z$  being solution values. This contradiction negates the assumption that $\gcd (z,p) = 1$. Conclude that $\gcd (z,p) = p$  because $p$  is prime.
\end{proof}


\begin{proof} \label{lem: five}
\underline{Lemma 5}: \\

\noindent Rewriting equation (\ref{eqn: se})  with the new notation and performing some algebra makes the equation $4x^{\circ}y^{\circ}z^{\circ}abcd = (x^{\circ}y^{\circ}a + x^{\circ}z^{\circ}b + y^{\circ}z^{\circ}c)p$.  Without loss of generality, suppose that a prime $q \neq p$  divides one of $x^{\circ}, y^{\circ}$,  and $z^{\circ}$, for example $q | x^{\circ}$.  The equation would imply that $q | y^{\circ}z^{\circ}c$; however, definitionally, $q \nmid y^{\circ}$ because $\gcd(x^{\circ},y^{\circ})=1$, $q \nmid z^{\circ}$ because $\gcd(x^{\circ},z^{\circ})=1$, and $q \nmid c$  because $\gcd(x^{\circ},c)=1$. Therefore $q \nmid x^{\circ}$.  The same will be true that primes $q \neq p$  has $p \nmid y^{\circ}$  and $p \nmid z^{\circ}$. 

For a Type I solution the prime $p \nmid x$  and $p \nmid y$.  This implies that $x^{\circ}=y^{\circ}=1$.  $p|z$, $p \nmid \gcd(x,z)$,  and $p \nmid \gcd(y,z)$  implies that $z^{\circ}=p$.  For a Type II solution the prime $p \nmid x$.  This implies that $x^{\circ}=1$.  Lemma~\ref{lem: c}  indicates that both $\gcd(y,p^{2}) \neq p^{2}$  and $\gcd(z,p^{2}) \neq p^{2}$.  Considering this along with $p|c$ implies that $p \nmid y^{\circ}$  and $p \nmid z^{\circ}$.  This means that $y^{\circ}=z^{\circ}=1$.
\end{proof}


\begin{proof} \label{thm: one}
\underline{Theorem 1}: \\

\noindent  Using the factorizations of $x,y,$  and $z$  with Lemma~\ref{lem: e}  makes $x=abd, y=acd, z=bcdp,$  and $p= (4abcd - a) \slash (b+c)$ for Type I solutions.  This makes

\begin{align*}
4xy - (x+y)p &= 4a^2bcd^2 - ad(4abcd - a) \\
&= a^{2}d.
\end{align*}

\ 

Let $c^{*} \in \mathbb{N}$  so that $c = c^{*}p$.  Using the factorizations of $x,y,$  and $z$  with Lemma~\ref{lem: e}  makes $x=abd, y=acd, z=bcd,$  and  $p= 4abd - (a+b) \slash c^{*}$  for Type II solutions.  This makes

\begin{align*}
4xy - (x+y)p &= p(4a^2bc^{*}d^2 - abd - ad(4abc^{*}d - (a+b))) \\
&= pa^{2}d.
\end{align*}

\ 

This implies that regardless of the type of solution, $4xy - (x+y)p = \gcd(y,p) a^{2}d$.  For a Type I solution $p(b+c) = a(4bcd - 1)$  and $\gcd(a,p)=1$  implies that $a|(b+c)$.  Suppose a prime $q| ((b+c) \slash a)$.   This implies that $q|(4bcd - 1)$.  If $q|bcd$, then $q|1$.  This implies that $\gcd(bcd, (b+c) \slash a) = 1$.  Conclude that

\begin{align*}
\gcd(y,p) \gcd ( xy, x+y ) &= \gcd ( a^{2}bcd^{2}, abd + acd ) \\
&= a^{2}d \cdot \gcd \left( bcd, \frac{b+c}{a} \right) \\
&= a^{2}d.
\end{align*}

\ 

For a Type II solution $pc^{*} = a(4bc^{*}d - 1) - b$  which implies that if $p|(4bc^{*}d-1)$,  then $p|b$.  $\gcd(p,b)=1$  implies that $\gcd(p,4bc^{*}d-1)=1$.  Also note that $\gcd(bc^{*}d,  4bc^{*}d -1) = 1$, so it is clear that $\gcd(bc^{*}dp, 4bc^{*}d-1)=1$.  Conclude that

\begin{align*}
\gcd(y,p) \gcd ( xy, x+y ) &= p \gcd ( a^{2}bc^{*}d^{2}p, abd + ac^{*}dp ) \\
&= p a^{2}d \cdot \gcd \left( bc^{*}dp, 4bc^{*}d-1 \right) \\
&= p a^{2}d.
\end{align*}

\ 

This implies that regardless of the type of solution, $\gcd(y,p) \gcd(xy, x+y) = \gcd(y,p) a^{2}d$.  Finally conclude that regardless of the type of solution $4xy-(x+y)p =  \gcd(y,p) \gcd(xy,x+y)$.  
\end{proof}


\begin{proof} \label{thm: two}
\underline{Theorem 2}: \\

\noindent Theorem~\ref{thm: a} implies that $4xy - (x+y)p = \gcd (xy, x+y)$  for any Type I solutions.  Dividing both sides by $4y-p$  shows that 

\begin{align*}
x - \frac{yp}{4y-p} &= \frac{\gcd ( xy, x+y)}{4y - p} \\
&< \frac{x+y}{4y-p} \\
&< \frac{2y}{2y + (2y-p)}.
\end{align*}

\ 

\noindent By definition $y > p \slash 2$,  so  $$x - \frac{yp}{4y-p} < 1.$$

\ 

\noindent Noting that $$\frac{\gcd(xy,x+y)}{4y-p} > 0$$

\ 

\noindent and $x \in \mathbb{N}$, conclude $$x -  \left\lceil \frac{yp}{4y-p} \right\rceil = 0.$$
\end{proof}


\begin{proof} \label{cor: one}
\underline{Corollary 1}: \\

\noindent Consider possible solutions with $x > \left\lceil  p \slash 2  \right\rceil$.  Recall from the proof of Lemma \ref{lem: c}, $2xp \slash (4x - p)$  decreases as the integer value $x \geq \left\lceil p \slash 4 \right\rceil$  increases.  This implies $2xp \slash (4x - p) < 2\left\lceil p \slash 2 \right\rceil p \slash (4 \left\lceil p \slash 2 \right\rceil - p) \leq p$  for integers $x > \left\lceil p \slash 2 \right\rceil$.  Lemma \ref{lem: a}  then suggests $y \leq \left\lfloor 2xp \slash (4x-p) \right\rfloor < p$.  Type II solutions require $p|y$.  It is impossible to have Type II solutions when $x > \left\lceil p \slash 2 \right\rceil$  because $y<p$. 

Solution values also require $x \leq y$.  If $x > \left\lceil  p \slash 2  \right\rceil$, then $y> \left\lceil p \slash 2 \right\rceil$.  Similarly, it can be inferred from the proof of Lemma \ref{lem: c}, $yp \slash (4y - p)$  decreases as the integer value $y \geq \left\lceil p \slash 4 \right\rceil$  increases.  This implies $yp \slash (4y - p) < \left\lceil p \slash 2 \right\rceil p \slash (4 \left\lceil p \slash 2 \right\rceil - p) \leq \left\lceil p \slash 2 \right\rceil$  for integers $y > \left\lceil p \slash 2 \right\rceil$.  Theorem \ref{thm: b} makes clear that for Type I solutions $x = \left\lceil yp \slash (4y-p) \right\rceil$, implying that $x \leq \left\lceil p \slash 2 \right\rceil$  for $y > \left\lceil p \slash 2 \right\rceil$.  Note that this cannot occur, because $x > \left\lceil p \slash 2 \right\rceil$.  Conclude that $x \leq \left\lceil p \slash 2 \right\rceil$, because there are are no Type I or Type II solutions for $x > \left\lceil p \slash 2 \right\rceil$.  Combining this with Lemma \ref{lem: a}  finishes the Corollary.\\
\end{proof}

\bibliographystyle{amsplain}

\end{document}